\documentclass[12pt]{article}
\usepackage[a4paper, total={7in, 10in}]{geometry}
\usepackage[parfill]{parskip}
\usepackage{physics, tensor, float, subcaption}
\usepackage{graphicx}
\graphicspath{ {Plots/} }
\usepackage{jhep-mod}
\usepackage{bm}
\usepackage{soul}
\usepackage{mathrsfs}
\usepackage[utf8]{inputenc}
\usepackage{enumerate}
\usepackage{bigints}
\usepackage{xcolor}
\usepackage{appendix}
\usepackage{graphicx}
\usepackage{float}
\usepackage{tikz}
\usepackage{setspace}
\usepackage{cancel}
\usepackage{doi}
\usepackage{amssymb,amsmath,amsthm}
\definecolor{purple}{rgb}{1,0,1}
\definecolor{lime}{HTML}{A6CE39} % needs xcolor
\newcommand{\blue}[1]{{\color{blue} #1}}
\newcommand{\B}[1]{{\bf\color{blue} #1}}

\definecolor{lime}{HTML}{A6CE39}
\newcommand{\orcidicon}{%
	\begin{tikzpicture}
	\draw[lime, fill=lime] (0,0) 
		circle [radius=0.16] 
		node[white] {{\fontfamily{qag}\selectfont \tiny ID}};
	\draw[white, fill=white] (-0.0625,0.095) 
		circle [radius=0.007];
	\end{tikzpicture}
	\hspace{-5mm}
}
\newcommand\orcidMatt{{\href{https://orcid.org/0000-0003-1088-6485}{\orcidicon}}}
\renewcommand{\O}{\mathcal{O}}
\def\floor{{\mathrm{floor}}}
\def\ceil{{\mathrm{ceil}}}

\def\Sierpinski{{Sierpi\'nski}}
\begin{document}
%========================================================
%========================================================
%========================================================

\title{\null\vspace{-50pt}
\leftline{\blue{ \Sierpinski's Hypothesis H$_1$}}}

%========================================================
%========================================================
%========================================================
\author{
\Large
Matt Visser\!\orcidMatt\!
}
%========================================================
%========================================================
%========================================================
%========================================================
\affiliation{School of Mathematics and Statistics, Victoria University of Wellington, \\
\null\qquad PO Box 600, Wellington 6140, New Zealand.}
%========================================================
%========================================================
\emailAdd{matt.visser@sms.vuw.ac.nz}
%========================================================
%========================================================
\def\theta{\vartheta}
\def\O{{\mathcal{O}}}
\def\Li{{\mathrm{Li}}}

\abstract{
\vspace{1em}

\Sierpinski's Hypothesis H$_1$, formulated in 1958, is the conjecture that (provided \break $n\geq 2$),  when the first $n^2$ counting numbers, $1, 2,3,\dots n^2$, are arranged in a square, then each row contains at least one  prime.
This conjecture is particularly interesting in that it subsumes and is stronger than both the Oppermann and Legrendre conjectures.
Herein, using the currently known maximal prime gaps,  I shall verify \Sierpinski's Hypothesis H$_1$ for (at least) the first $n \leq \hbox{10 070 368 414} \gtrsim 10 \hbox{ billion}$ of these \Sierpinski{} matrices. 
I~shall also demonstrate some partial but more general results. For example: Even for arbitrary $n\geq \hbox{10 070 368 414}$  at least one quarter of the rows of the $n$th \Sierpinski{} matrix contain at least one prime. 
Furthermore, even for arbitrary $n\geq \hbox{10 070 368 414}$ at least the first $\hbox{141 618}$ rows of the $n$th \Sierpinski{} matrix always contain at least one prime. 
These and related results are obtained largely by using the locations and values of the known maximal prime gaps, the pigeonhole principle, and some recent bounds on the first Chebyshev function.

\bigskip

\bigskip
\noindent
{\sc Date:} \\
Saturday 27 December 2025; Updated 20 July 2026; \LaTeX-ed \today


\bigskip
\noindent{\sc Mathematics Subject Classification 2020}: \\
\null \qquad 11A41: Primes; 11N05: Distribution of primes.

\bigskip
\noindent{\sc Changes}: Various minor typos fixed; Presentation improved; \\
Updated to reflect the discovery of the 85th maximal prime gap. 

}

%========================================================
\maketitle
%========================================================
\def\tr{{\mathrm{tr}}}
\def\diag{{\mathrm{diag}}}
\def\cof{{\mathrm{cof}}}
\def\pdet{{\mathrm{pdet}}}
\def\d{{\mathrm{d}}}
\parindent0pt
\parskip7pt
\def\Kerr{{\scriptscriptstyle{\mathrm{Kerr}}}}
\def\eos{{\scriptscriptstyle{\mathrm{eos}}}}

\def\Z{{\mathbb{Z}}}

\clearpage

\null
\vspace{-60pt}
%================================================
\section{Introduction}
%================================================

\Sierpinski's Hypothesis H$_1$ was introduced in the classic 1958 article 
``Sur certaines hypoth\`eses concernant les nombres premiers'',
(On certain hypotheses concerning the prime numbers), authored by \Sierpinski{}  and Schinzel~\cite{S_and_S}.  
See also reference~\cite{Remarks}.

{\bf Hypothesis H$_1$:} \emph{ Hypothesis H$_1$ is the conjecture that (provided  $n\geq 2$) when the first $n^2$ counting numbers, $1, 2,3,\dots n^2$, are arranged in a square, then each row contains at least one  prime.}
\hfill $\Box$.

Specifically let the matrix $S_n$ be given by
\begin{equation}
[S_n]_{ij} = n (n-i) + j; \qquad i,j\in\{1,\dots,n\}; \qquad n\geq 2. 
\end{equation}
That is
\begin{equation}
S_n = 
\begin{bmatrix} 
(n-1)n+1 & (n-1)n+2 & \cdots &  n^2-1 & n^2\\
(n-2)n+1 & (n-2)n+2 & \cdots &  (n-1)n-1 & (n-1)n\\
\vdots &\vdots &\ddots &\vdots &\vdots\\
n+1 & n+2 & \cdots &  2n-1 & 2n\\
1 & 2 & \cdots & n-1 & n
 \end{bmatrix}
\end{equation}
We can rewrite this as
\begin{equation}
S_n = 
\begin{bmatrix} 
n^2-n+1 & n^2-n+2 & \cdots &  n^2-1 & n^2\\
(n-1)^2 & (n-1)^2 +1 & \cdots &  n^2-n-1 & n^2-n\\
\vdots &\vdots &\ddots &\vdots &\vdots\\
n+1 & n+2 & \cdots &  2n-1 & 2n\\
1 & 2 & \cdots & n-1 & n\\
 \end{bmatrix}
\end{equation}
\enlargethispage{20pt}
\begin{itemize}
\itemsep-5pt
\item 
That the top row contains at least one prime is equivalent to the first half of Oppermann's conjecture.
\item
That the second to top row contains at least one prime is equivalent to the second half of Oppermann's conjecture.
\item
That the top two rows between them contain at least one  prime is equivalent to Legendre's conjecture.
\item
That the bottom row contains a prime is trivial.
\item
That the second row from the bottom contains a prime is guaranteed by the Bertrand--Chebyshev theorem.
\item
 That for $n \geq$ 9 each of the lowest 9 rows contains at least one prime number can easily be deduced from the result that for 
 $n \geq 5$ there is at least one prime number in the interval $(8n,9n)$.
 (This is a variation on  Breusch's 1932 result~\cite{Breusch} that for $x \geq 48$ there is at least one prime number between $x$ and $\frac98 x$.) Below we shall see how to replace the number 9 above by \hbox{141 618}.
\item
That Hypothesis H$_1$ certainly holds for $n\in[2,\hbox{ 4 505}]$ was pointed out by Schinzel in 1961. See reference~\cite{Remarks}.
\end{itemize}

\enlargethispage{20pt}
Below I shall, to the extent possible, derive expanded explicit regions of validity for Hypothesis H$_1$. The relevant tools to be used involve the locations and values of the known maximal prime gaps, the pigeonhole principle, and some recent bounds on the first Chebyshev function.

%\clearpage

\enlargethispage{20pt}
To set the stage, here are a few specific examples of \Sierpinski{} matrices with the primes emphasized in boldface blue:
\begin{equation}
S_2 = \begin{bmatrix}  
\B{3} &4\\1 & \B{2}  
\end{bmatrix};
\qquad
S_3 = \begin{bmatrix} 
\B{7}&8&9\\ 4&\B{5}&6\\ 1&\B{2}&\B{3} 
\end{bmatrix};
\qquad
S_4 = \begin{bmatrix} 
\B{13}&14&15&16\\ 9&10&\B{11}&12\\ \B{5}&6&\B{7}&8\\1&\B{2}&\B{3}&4
\end{bmatrix};
\end{equation}
and
\begin{equation}
S_5 = \begin{bmatrix} 
21&22&\B{23}&24&25\\
16&\B{17}&18&\B{19}&20\\
\B{11}&12&\B{13}&14&15\\
6&\B{7}&8&9&10\\
1&\B{2}&\B{3}&4&\B{5}
\end{bmatrix}.
\end{equation}
You can see that the number of primes in each row fluctuates quite markedly.\\
With a bit more (elementary) work one sees that 
%{\small
\begin{equation}
S_{10} = \begin{bmatrix} 
91&92&93&94&95&96&\B{97}&98&99&100\\
81&82&\B{83}&84&85&86&{87}&88&89&90\\
\B{71}&72&\B{73}&74&75&76&77&78&\B{79}&80\\
\B{61}&62&63&64&65&66&\B{67}&68&69&70\\
51&52&\B{53}&54&55&56&57&58&\B{59}&60\\
\B{41}&42&\B{43}&44&45&46&\B{47}&48&49&50\\
\B{31}&32&33&34&35&36&\B{37}&38&39&40\\
21&22&\B{23}&24&25&26&27&28&\B{29}&30\\
\B{11}&12&\B{13}&14&15&16&\B{17}&18&\B{19}&20\\
1&\B{2}&\B{3}&4&\B{5}&6&\B{7}&8&9&10
\end{bmatrix}
\end{equation}
%}
and 
{\small
\begin{equation}
S_{15} = \left[\begin{array}{ccccccccccccccc}
\B{211}&212&213&214&215&216&217&218&219&220&221&222&\B{223}&224&225\\
196&\B{197}&198&\B{199}&200&201&202&203&204&205&206&207&208&209&210\\
\B{181}&182&183&184&185&186&187&188&189&190&\B{191}&192&\B{193}&194&195\\
166&\B{167}&168&169&170&171&172&\B{173}&174&175&176&177&178&\B{179}&180\\
\B{151}&152&153&154&155&156&\B{157}&158&159&160&161&162&\B{163}&164&165\\
136&\B{137}&138&\B{139}&140&141&142&143&144&145&146&147&148&\B{149}&150\\
121&122&123&124&125&126&\B{127}&128&129&130&\B{131}&132&133&134&135\\
106&\B{107}&108&\B{109}&110&111&112&\B{113}&114&115&116&117&118&119&120\\
91&92&93&94&95&96&\B{97}&98&99&100&\B{101}&102&\B{103}&104&105\\
76&77&78&\B{79}&80&81&82&\B{83}&84&85&86&87&88&\B{89}&90\\
\B{61}&62&63&64&65&66&\B{67}&68&69&70&\B{71}&72&\B{73}&74&75\\
46&\B{47}&48&49&50&51&52&\B{53}&54&55&56&57&58&\B{59}&60\\
\B{31}&32&33&34&35&36&\B{37}&38&39&40&\B{41}&42&\B{43}&44&45\\
16&\B{17}&18&\B{19}&20&21&22&\B{23}&24&25&26&27&28&\B{29}&30\\
1&\B2&\B3&4&\B5&6&\B7&8&9&10&\B{11}&12&\B{13}&14&15\\
\end{array}\right]
\end{equation}
}
You can see that the number of primes in each row fluctuates quite markedly.

Note that on average the rows of the \Sierpinski{} matrix $S_n$ will contain approximately $\pi(n^2)/n \sim \frac12\pi(n)$ primes, so it will be the fluctuations (not the averages) that determine the validity or otherwise of Hypothesis H$_1$. 

\bigskip
\hrule\hrule\hrule

%================================================
\section{Two stronger hypotheses}
%================================================
\enlargethispage{40pt}
It is relatively straightforward to develop conjectures stronger than \Sierpinski's H$_1$. For instance Schinzel has suggested~\cite{Remarks}:
\begin{itemize}
\item 
{\bf C1:} The conjecture (only slightly stronger than Hypothesis H$_1$)
\begin{equation}
\forall x\geq 117,\quad \exists\; p\in(x,x+\sqrt{x}) 
\end{equation}
implies \Sierpinski's H$_1$.
\item 
{\bf C2:}  The conjecture (considerably stronger than Hypothesis H$_1$)
\begin{equation}
\forall x\geq 8,\quad \exists\; p \in(x,x+\ln^3 x) 
\end{equation}
implies \Sierpinski's H$_1$.
\end{itemize}
Known ranges of validity for these stronger conjectures were used by Schinzel in 1961~\cite{Remarks} to argue that H$_1$ is valid at least for $n\in[2,\;\hbox{4 505}]$. While we will use this specific validity interval in one of the  arguments below, we will not directly make use of these stronger conjectures.

%\clearpage
\bigskip
\hrule\hrule\hrule
%================================================
\section
{The situation for  $n\leq \hbox{10 070 368 414}$}
%================================================
For the rather sizeable region $n\in[2,\;\hbox{10 070 368 414}] \supset [2, 10^{10}]$ we can extract full information regarding Hypothesis H$_1$.

%================================================
\subsection{Unconditional validity for H$_1$. }
%================================================

{\bf Theorem 1:}
\Sierpinski's Hypothesis H$_1$ holds  for all \Sierpinski{} matrices $S_n$ from $n=2$ at least up to $n=\hbox{10 070 368 414}\gtrsim 10^{10}$. 

{\bf Proof technique:} Use of the locations and sizes of the known maximal prime gaps to cover the region $n\in[\hbox{1 572},\; \hbox{10 070 368 414}]$, and and Schinzel's result in the overlapping region $n\in[2,\;\hbox{4 505}]$.

{\bf Proof:}
The current (as of July 2026)  largest known maximal prime gap is the 85th maximal prime gap, located at $p_{85}^*= \hbox{101 412 319 996 363 309 069} \approx 10^{20}  \gtrsim 2^{66}$, 
with $n^*_{85}= \pi(p^*_{85})=2 251 483 061 895 611 799
\approx 2.25\times 10^{18}\gtrsim 2^{60}$, and gap size $g_{85}^*=\hbox{1 854}$. 
See references~\cite{prime-gaps, tables, primecount}. 

Define the number $N_{85} = \floor[\sqrt{p_{85}^*}\,] =\hbox{10 070 368 414} \gtrsim 10^{10}$, and consider the $N_{85}\times N_{85}$ \Sierpinski{} matrix $S_{N_{83}}$ which contains $\hbox{101 412 319 993 688 875 396} < p^*_{85}$ elements.

Thence the maximum gap between primes within the matrix $S_{N_{85}}$ will be governed by the  84th maximal prime gap, which is located at $p_{84}^*= \hbox{
68 068 810 283 234 182 907} \approx 6.8\times 10^{19}  \gtrsim 2^{65}$, with $n^*_{84}= \pi(p^*_{84})=\hbox{1 524 717 378 371 224 128}\approx 1.5\times 10^{18}\gtrsim 2^{60}$, and gap size $g_{84}^*=\hbox{1 724}$. 
See references~\cite{prime-gaps, tables, primecount}.

Specifically, by the pigeonhole principle, each specific row in this specific \Sierpinski{} matrix $S_{N_{85}}$ must contain, at the very least,
$\floor[N_{85}/g_{84}^*]=\floor[N_{85}/(\hbox{1 724})] = \hbox{5 841 280} > 5.8 \hbox{ million}$ primes.

Furthermore, for $n \in [g_{84}^*, N_{85}] = [\hbox{1 724}, N_{85}]$, each row in the \Sierpinski{} matrix $S_n$ must contain at least $\floor[n/(\hbox{1 724})] $ primes, certainly more than 1 prime.

This verifies \Sierpinski's hypothesis H$_1$ for all $n \in [\hbox{1 724}, N_{85}]$.

But this overlaps with the region of validity $n\in[2,\hbox{ 4 505}]$ already established by
Schinzel in 1961~\cite{Remarks}. 

Thus this now verifies \Sierpinski's hypothesis H$_1$ for all $n \in [2, N_{85}]$.

That is, the first 10 billion \Sierpinski{} matrices (at least) satisfy Hypothesis H$_1$.
\rightline{$\Box$.}

\enlargethispage{20pt}
{\bf Comment:} Note that the way the proof is set up it proves something slightly stronger than H$_1$. 
Not only do all of the individual rows in the \Sierpinski{} matrix contain primes, but any sequence of $n$ contiguous numbers in  the \Sierpinski{} matrix must contain at least one prime. 

\bigskip
\hrule\hrule\hrule
%%%%%%%%%%%%%%%%%%%%%%%%%%%%%%%%%%%%%%%%%%%%%%
\subsection{Alternate proof strategy: Descent}

Instead of appealing to Schinzel's result one could iterate the preceding argument for $n \in [g_{84}^*, N_{85}] = [\hbox{1 724}, N_{85}]$ downwards as far as possible. Doing so provides some extra information not captured by the previous argument.

{\bf Step 2:}
Take $n=g_{84}^*-1=\hbox{ 1 723}$ and consider $S_{\scriptsize{\hbox{1 723}}}$. Note $(\hbox{1 723})^2 = \hbox{2 968 729}$
which lies between the 21st and 22nd maximal prime gaps.
Indeed the 21st maximal prime gap is located at $p_{21}^*= \hbox{
2 010 733}$, with $n^*_{21}= \pi(p^*_{21})=\hbox{149 689}$, 
and gap size $g_{21}^*=148$. Similarly
 the 22nd maximal prime gap is located at $p_{22}^*= \hbox{4 652 353}$, with $n^*_{22}= \pi(p^*_{22})=\hbox{325 852}$, and gap size $g_{22}^*=154$. See references~\cite{prime-gaps, tables, primecount}.

Thence, by the pigeonhole principle, each row in the \Sierpinski{} matrix $S_{\scriptsize{\hbox{1 723}}}$ must contain 
at least
$\floor[\hbox{1 723}/148] = 11$ primes.
Furthermore, for $n \in [148, \hbox{1 723}]$ each row in the \Sierpinski{} matrix $S_n$ must contain at least $\floor[n/148] \geq 1$ prime.

This now verifies \Sierpinski's hypothesis H$_1$ for all $n \in [148, N_{85}]$.

%\clearpage
%%%%%%%%%%%%%%%%%%%%%%%%%
%%% 
{\bf Step 3:}
Now take $n=147$ and consider $S_{147}$. Note $(147)^2 = \hbox{21 609}$
which lies between the 13th and 14th maximal prime gaps.
Indeed the 13th maximal prime gap is located at $p_{13}^*= \hbox{
19 609}$, with $n^*_{13}= \pi(p^*_{13})=\hbox{254}$, 
and gap size $g_{21}^*=52$. 
Similarly
 the 14th maximal prime gap is located at $p_{14}^*= \hbox{31 397}$, with $n^*_{22}= \pi(p^*_{22})=\hbox{3 385}$, and gap size $g_{14}^*=72$.

Thence, by the pigeonhole principle, each row in the \Sierpinski{} matrix $S_{147}$ must contain 
at least
$\floor[147/52] = 2$ primes.
Furthermore, for $n \in [52,147]$ each row in the \Sierpinski{} matrix $S_n$ must contain at least $\floor[n/52] \geq 1$ prime.

This now verifies \Sierpinski's hypothesis H$_1$ for all $n \in [52, N_{85}]$.

%%%%%%%%%%%%%%%%%%%%%%%%%

{\bf Step 4:} 
Now take $n=51$ and consider $S_{51}$. Note $(51)^2 = \hbox{2601}$
which lies between the 10th and 11th maximal prime gaps.
Indeed the 10th maximal prime gap is located at $p_{10}^*= \hbox{
1 327}$, with $n^*_{10}= \pi(p^*_{10})=\hbox{217}$, 
and gap size $g_{10}^*=34$. 
Similarly
 the 11th maximal prime gap is located at $p_{11}^*= \hbox{9551}$, with $n^*_{11}= \pi(p^*_{11})=\hbox{1 183}$, and gap size $g_{11}^*=36$.
 See references~\cite{prime-gaps, tables, primecount}.

Thence, by the pigeonhole principle, each row in the \Sierpinski{} matrix $S_{51}$ must contain 
at least
$\floor[51/34] = 1$ prime.
Furthermore, for $n \in [34,51]$ each row in the \Sierpinski{} matrix $S_n$ must contain at least $\floor[n/34] = 1$ prime.

This now verifies \Sierpinski's hypothesis H$_1$ for all $n \in [34, N_{85}]$.

%%%%%%%%%%%%%%%%%%%%%%%%%

{\bf Step 5:} 
Now take $n=33$ and consider $S_{33}$. Note $(33)^2 = \hbox{1089}$
which lies between the 8th and 9th maximal prime gaps.
Indeed the 8th maximal prime gap is located at $p_{8}^*= \hbox{
887}$, with $n^*_{8}= \pi(p^*_{8})=\hbox{154}$, 
and gap size $g_{8}^*=20$. 
Similarly
 the 9th maximal prime gap is located at $p_{9}^*= \hbox{1129}$, with $n^*_{9}= \pi(p^*_{9})=\hbox{189}$, and gap size $g_{9}^*=22$.
 See references~\cite{prime-gaps, tables, primecount}.

Thence each row in the \Sierpinski{} matrix $S_{33}$ must contain 
at least
$\floor[33/20] = 1$ prime.
Furthermore, for $n \in [20,33]$ each row in the \Sierpinski{} matrix $S_n$ must contain at least $\floor[n/20] = 1$ prime.

This now verifies \Sierpinski's hypothesis H$_1$ for all $n \in [20, N_{85}]$.

%\clearpage
%%%%%%%%%%%%%%%%%%%%%%%%%
{\bf Step 6:} 
Now take $n=19$ and consider $S_{19}$. Note $(19)^2 = \hbox{361}$
which lies between the 6th and 7th maximal prime gaps.
Indeed the 6th maximal prime gap is located at $p_{6}^*= \hbox{
113}$, with $n^*_{6}= \pi(p^*_{6})=\hbox{30}$, 
and gap size $g_{6}^*=14$. 
Similarly
 the 7th maximal prime gap is located at $p_{7}^*= \hbox{523}$, with $n^*_{7}= \pi(p^*_{7})=\hbox{99}$, and gap size $g_{7}^*=18$.
 See references~\cite{prime-gaps, tables, primecount}.

Thence, by the pigeonhole principle, each row in the \Sierpinski{} matrix $S_{19}$ must contain 
at least
$\floor[19/14] = 1$ prime.
Furthermore, for $n \in [14,19]$ each row in the \Sierpinski{} matrix $S_n$ must contain at least $\floor[n/14] \geq 1$ prime.

This now verifies \Sierpinski's hypothesis H$_1$ for all $n \in [14, N_{85}]$.

%%%%%%%%%%%%%%%%%%%%%%%%%
{\bf Step 7:} 
Now take $n=13$ and consider $S_{13}$. Note $(13)^2 = \hbox{167}$
which again lies between the 6th and 7th maximal prime gaps.

Thence the arguments we have given above now yield no further information: they merely indicate that each row in the \Sierpinski{} matrix $S_{13}$ must contain 
at least
$\floor[13/g_6^*] = \floor[13/14] = 0$ primes.

Fortunately $13\times 13$ matrices (and smaller) are not too difficult to deal with ``by hand''.  We have  already explicitly checked $S_2$, $S_3$, $S_4$, $S_5$,  $S_{10}$,  and $S_{15}$ in the Introduction. 
For the remaining cases, since the first two rows of the matrix are trivial, a  simple diagnostic is to check whether
\begin{equation}
Q_n = \prod_{j=2}^{n-1} \left\{ \pi([j+1]n)-\pi(jn) \right\} > 0.
\end{equation}
Checking these inequalities for $n\in[2,13]$  now verifies \Sierpinski's Hypothesis H$_1$ for all $n \in [2, N_{85}]$, where $N_{85} = \floor[\sqrt{p_{85}^*}\,] =\hbox{10 070 368 414}$. 

That is, the first 10 billion \Sierpinski{} matrices (at least) satisfy Hypothesis H$_1$.
\rightline{$\Box$.}

\bigskip
\hrule\hrule\hrule

\clearpage
%================================================
\subsection{Lower bounds on the number of primes in each row. }
%================================================

\enlargethispage{40pt}
A side effect of the discussion above is that for any $n\leq N_{85} \approx 10^{10}$ we have rather good control on the lowest possible  number of primes in each row.

{\bf Theorem 2:} For selected values of $n$, lower bounds on the  minimum number of primes in each row of the \Sierpinski{} matrix $S_n$ 
are reported in Table I below.

{\bf Strategy:} For any $n\in[2,N_{85}]$ first calculate $n^2$ and then, using tables of known maximal prime gaps~\cite{prime-gaps, tables},
 find the location of the largest maximal prime gap with starting point below $n^2$. Define the index function $i(n) = \max\{i: p_i^*< n^2\}$. Then every prime gap below $n^2$ is certainly less than or equal to $g_{i(n)}^*$ in width, and therefore any $n$ contiguous numbers in the $n$th \Sierpinski{} matrix $S_n$ must contain at least $\#(n) = \floor [n/g_{i(n)}^*]$ primes. 
 
In particular, any specific row of the $n$th \Sierpinski{} matrix $S_n$ must contain at least $\#(n) = \floor [n/g_{i(n)}^*]$ primes. 
Selected bounds of this type are reported in Table I below.

\vspace{-10pt}
\begin{table}[!h]
\caption{Minimum number of primes in each row of $S_n$.}
\begin{center}
\begin{tabular}{||c|c|c|c|c||}
\hline\hline
$n$ & $n^2$ & $i(n)$ & $g_{i(n)}^*$ & $\#(n)$ \\
\hline\hline
$10^{10}$          & $10^{20}$                 &84 & 1 724 & 5 800 464\\
$6\times 10^9$ & $3.6 \times 10^{19}$ &83 & 1 676 & 3 579 952\\
$4\times 10^9$ & $1.6 \times 10^{19}$ &78 & 1 526 & 2 621 231\\
$3\times 10^9$ & $9 \times 10^{18}$    &77 & 1 510 & 1 986 754\\
$2\times 10^9$ & $4 \times 10^{18}$    &75 & 1 476 & 1 355 013\\
$10^9$             &  $ 10^{18} $               & 74  & 1 442 & 693 481\\
\hline
$10^8$ & $ 10^{16} $               &64& 1 132& 88 339\\ % 108 255\\ ?
$10^7$ & $ 10^{14} $               &59 & 804& 12 437\\
$10^6$ & $ 10^{12} $               &50& 540& 1 851 \\ %1 821 \\
$10^5$ & $ 10^{10} $               &35& 354& 282\\
$10^4$ & $ 10^{8} $                 &25& 220& 45 \\ %28\\
$10^3$ & $ 10^{6} $                 &18& 114& 8\\
$10^2 $ & $ 10^{4} $                &11&36& 2\\
\hline
$[23,25]$ &[529,625] &7 & 18& 1\\
$[14,22]$ &[196,484] &6 & 14& 1\\
\hline
\hline
[11,13]  & [121,169]&  6& 14& 0\\
10        &   100      & 5& 8& 1\\{}
[6,9]     & [36,81]   & 4&  6& 1\\
5          & 25          & 4&  6& 0\\{}
4          & 16          & 3& 4& 1\\{}
3          & 9            & 3& 4& 0\\{}
2          &4             & 2& 2& 1\\
%$2$ & 4& 4& 0\\
\hline\hline
\end{tabular}
\end{center}
\label{T:table1}
\end{table}%

Note that, (in analogy with what we already saw happening for the descent version of Theorem 1), for $n\leq 13$ this algorithm often fails to provide useful information, and one must resort to explicit calculation.

\bigskip
\hrule\hrule\hrule

\clearpage
%================================================
%================================================
\section{The situation for  $n\geq \hbox{10 070 368 414}$ }
%================================================
What if anything can we say for $n>N_{85}$? We shall develop a number of partial but nevertheless interesting results. 
%================================================
\subsection{Fractional results (for possibly non-contiguous rows).}
%================================================
{\bf Theorem 3:}
For arbitrary $n$, at least one quarter of the rows of the \Sierpinski{} matrix $S_n$ contain at least one prime. 

{\bf Proof:} 
For $n\leq N_{85}$ a stronger result has already been established, \emph{all} the rows contain primes. For $n>N_{85}$ we 
 use the classic result that $\pi(x) > x/\ln x$ for reals $x\geq 17$~\cite{Bounds}, which on integers implies  $\pi(n) > n/\ln n$ for integers $n\geq 11$. So certainly $\pi(n^2) > {n^2/\ln(n^2)} = {n^2/(2\ln n)}$ over the region of interest. Furthermore, from the Montgomery--Vaughn bound
$\pi(x+y) < \pi(x)+ {2y\over\ln y}$ obtained from the large sieve~\cite{not-HL2},  we have $\pi([k+1]n)-\pi(kn) < {2n/\ln n}$.

That is, $n$th \Sierpinski{} matrix $S_n$ contains at least $ {n^2/(2\ln n)}$ primes, while each individual row in $S_n$ contains at most ${2n/\ln n}$ primes.
So quantitatively, by the pigeonhole principle,  at least
\begin{equation}
{ {n^2/(2\ln n)} \over {2n/\ln n}} = \frac{n}4
\end{equation}
of the rows of the $n$th \Sierpinski{} matrix must contain primes. That is,  working on a fractional basis, at least one quarter of the rows of any \Sierpinski{} matrix have to contain primes. \hfill{$\Box$.}

Unfortunately the current argument gives no information on \emph{where} these prime-containing rows might be located.
To slightly improve on this let us consider the bottom $m\leq n$ rows of the \Sierpinski{} matrix.

{\bf Theorem 4:}
For arbitrary $n$, at least a fraction $ \frac12 \, {\ln n\over\ln n+\ln m}$ of the bottom $m$ rows of the \Sierpinski{} matrix $S_n$ have to contain primes. 

{\bf Proof:} 
For $n\leq N_{85}$ a stronger result has already been established, \emph{all} the rows contain primes. For $n>N_{85}$ and $m\in[1,n]$ we  certainly have $\pi(m n) > {m n/\ln(mn)}$.
That is, the bottom $m$ rows of $S_n$ contain at least $ {m n/\ln(mn)}$ primes. 

\enlargethispage{40pt}
By the pigeonhole principle, quantitatively at least
\begin{equation}
{ {mn/\ln(mn)} \over {2n/\ln n}} = \frac{m}2 \; {\ln n\over\ln n+\ln m}
\end{equation}
of these bottom $m$ rows of the $n$th \Sierpinski{} matrix have to contain primes. \\
That is, on a fractional basis, at least 
\begin{equation}
 \frac12 \; {\ln n\over\ln n+\ln m}
\end{equation}
 of the bottom $m$ rows of the $n$th \Sierpinski{} matrix have to contain primes. 
\hfill $\Box$.

Note that as $m\to n$ Theorem 3 reduces to Theorem 2.

Unfortunately the current argument still gives no precise information on \emph{where} these prime-containing rows might be located.
We shall somewhat improve on this in the discussion below.

\bigskip
\hrule\hrule\hrule

%================================================
\subsection{Results for contiguous rows.}
%================================================

{\bf Theorem 5:}
Let $n>N_{85}$, but assume $n< \frac12(N_{85})^2$. Then at the very least all of the $\floor[(N_{85})^2/n]$ bottom rows of the \Sierpinski{} matrix $S_n$ contain primes.

{\bf Proof:}
We have already seen that the maximal prime gap in the interval $[1,(N_{85})^2]$ is 1 724.
Consequently any $n$ consecutive numbers in the interval $[1,(N_{85})^2]$ must contain at least $\floor(n/\hbox{1 724})$ primes. So the rows in $S_n$ up to but not including the row containing the number $ (N_{85})^2$ must contain at least $\floor(n/\hbox{1 724})$ primes. That is to say, the bottom 
$\floor[(N_{85})^2/n]$ rows of the of the \Sierpinski{} matrix $S_n$ must each contain at least $\floor(n/\hbox{1 724})$ primes. Note that this argument ceases to provide useful information once $n\geq \frac12(N_{85})^2$.
\hfill\hfill $\Box$.

{\bf Lemma:}
\Sierpinski's hypothesis H$_1$ is equivalent to the doubly infinite chain of strict inequalities
\begin{equation}
\pi([k+1]n)> \pi(kn) ; \qquad (\forall\; k \geq 1 \hbox{ and } n\geq k+1).
\end{equation}

{\bf Proof:} By inspection. \hfill\hfill $\Box$.

{\bf Theorem 6:}
 If $n$ is arbitrary, and $k\leq N_{85}-3 = \hbox{10 070 368 411} $, and the $k$th row of $S_n$ contains a prime, then all lower rows of $S_n$ also contain at least one prime. 
 
 {\bf Remark:} Of course this result is interesting only for $n>N_{85}$, and  in order for this result to be useful you still somehow have to verify that $k$th row of $S_n$ contains a prime.

{\bf Proof:}
Suppose H$_1$ holds for the $k$th row of $S_n$. \\
That is, suppose for some fixed $k\geq 1$ and $n\geq k+1$ we have
\begin{equation}
\pi([k+1]n)> \pi(kn) ; \qquad 
(\forall\; k \geq 1 \hbox{ and } n\geq k+1).
\end{equation}
Then setting $\tilde n = kn$ we see
\begin{equation}
\pi\left(\tilde n\left[1+ \frac1k\right]\right)> \pi(\tilde n) ; \qquad 
(k \geq 1 \hbox{ and } \tilde n\geq k[k+1] ).
\end{equation}
Then automatically we have the weaker result
\begin{equation}
\pi\left(\tilde n\left[1+ \frac1{k-1}\right]\right)> \pi(\tilde n) ; \qquad 
(k \geq 1 \hbox{ and } \tilde n\geq k[k+1] ).
\end{equation}
Thence
\begin{equation}
\pi\left(\tilde n\left[\frac{k}{k-1}\right]\right)> \pi(\tilde n) ; \qquad 
(k \geq 1 \hbox{ and } \tilde n\geq k[k+1] ).
\end{equation}
Then setting $\tilde n= \bar n [k-1]$
\begin{equation}
\pi\left(\bar{n} k \right)> \pi(\bar{n} [k-1]) ; \qquad 
(k \geq 2 \hbox{ and } \bar{n}\geq {k[k+1]/[k-1]} ).
\end{equation}
Relabelling $\bar n \to n$ this is \emph{almost} the $(k-1)$th row of the Sierpinski{} matrix $S_n$ which would be
\begin{equation}
\pi\left({n} k \right)> \pi({n} [k-1]) ; \qquad 
({n}\geq {k} ).
\end{equation}
Now
\begin{equation}
k < {k[k+1]\over [k-1]} \leq k+3 \qquad (k\geq 3)
\end{equation}
So in order for the $k$th row to imply the $(k-1)$th row  we would ``merely'' need to check the additional inequalities
\begin{equation}
\pi\left({n} k \right)> \pi({n} [k-1]) ; \qquad 
( {n}\in [k,k+3] ).
\end{equation}
This means we need to check the four inequalities
\begin{equation}
\pi\left(k^2 \right)> \pi(k [k-1]) ; \qquad 
\pi\left([k+1]k \right)> \pi([k+1] [k-1]) ; 
\end{equation} 
\begin{equation}
\pi\left([k+2]k \right)> \pi([k+2] [k-1]) ; \qquad
\pi\left([k+3]k \right)> \pi([k+3] [k-1]) ,
\end{equation}
But these inequalities correspond respectively to the top row of $S_k$, the second to top row of $S_{k+1}$, the third to top row of $S_{k+2}$, and the fourth to top row of $S_{k+3}$.

So as long as $k+3 \leq N_{85}$ we know all four of these inequalities are satisfied.

Consequently, let $n$ be arbitrary. (The case $n\leq N_{85}$ has already been fully dealt with; it is  the case $n>N_{85}$ that is of interest.) Then: If $k\leq N_{85}-3$ and the $k$th row of $S_n$ satisfies H$_1$, then the $(k-1)$th row of $S_n$ also satisfies H$_1$.

Consequently, if $k\leq N_{85}-3$ and the $k$th row of $S_n$ contains a prime,
then all lower rows of $S_n$ also contain a prime. (Of course in order for this to be useful you somehow have to verify that $k$th row of $S_n$ contains a prime.) \hfill $\Box$.

%\clearpage
\bigskip
\hrule\hrule\hrule

\clearpage
%================================================
\subsection{Bounds derived from the first Chebyshev function.}
%================================================

{\bf Theorem 7:} For all $n\geq \hbox{141 618}$ the $k$th row of the \Sierpinski{} matrix $S_n$ certainly contains a prime  for $1 \leq k \leq \hbox{141 618}$.

{\bf Proof:}
There are very many well-known  bounds of the form
\begin{equation}
|\theta(x)-x| < {\eta_m \, x \over \ln^m x } \qquad (x>x_m)
\end{equation}
We will avoid any specific choice of parameters $\{m,\eta_m,x_m\}$ for now.

We first note that we would want to apply these bounds in the form
\begin{equation}
|\theta(kn)-kn| < {\eta_m \; kn \over \ln^m(kn) } \qquad (kn>x_m).
\end{equation}
Now observe that
\begin{equation}
\theta([k+1]n)-\theta(kn) 
> n - {\eta_m [k+1]n \over\ln^m([k+1]n)} - {\eta_m kn \over \ln^m(kn)}
>  n - {\eta_m [k+1]n \over\ln^m(kn)} - {\eta_m kn \over \ln^m(kn)}
\end{equation}
\begin{equation}
= n - {\eta _m[2k+1] n\over \ln^m(kn)}= n \left( 1- {\eta _m[2k+1] \over \ln^m(kn)}\right)
\end{equation}
This is positive, thereby implying the existence of at least one prime in the interval $[nk, n[k+1]]$, whenever
\begin{equation}
\ln^m(kn) > \eta _m[2k+1] \qquad\hbox{subject to}\qquad k>{x_m\over n}.
\end{equation}
We already have full information for $n\leq N_{85}$, so let us consider $n>N_{85}$. Then certainly the $k$th row contains a prime provided
\begin{equation}
\ln^m(kN_{85}) > \eta _m[2k+1] \qquad\hbox{and}\qquad k>{x_m\over N_{85}}
\label{E:ineq}
\end{equation}
Some experimentation is now required (using the tables in Dusart~\cite{Dusart} and Axler~\cite{Axler}) to find the most useful combination of parameters $\{m, \eta_m, x_m\}$. 

\enlargethispage{40pt}
\begin{itemize}
\item 
One useful choice is to take  $m=1$, $\eta_1 = 0.001$, $x_1= \hbox{908 994 923}= p_{46\;445\;479}$ [Dusart~\cite{Dusart}]. The inequality  (\ref{E:ineq}) is then easily checked to be  satisfied on the interval $k\in[1,\; 16 367]$.
That is, for all $n>N_{83}$ the $k$th row of the \Sierpinski{} matrix $S_n$ certainly contains a prime  for $1 \leq k \leq \hbox{16 367}$.

 \item 
 Another useful choice is to take  $m=2$, $\eta_2 = 0.01$, and $x_2= \hbox{7 713 133 853}= p_{355\;191\;539}$ [Dusart~\cite{Dusart}].
The inequality (\ref{E:ineq}) is valid only for $k\geq \ceil[x_2/N_{85}]=1$, and is  then
easily checked to be satisfied on the interval $k\in[1,\; \hbox{57 790}]$. 
%Of course the interval $k\in[1,9]$ was already dealt with in \Sierpinski's original discussion~\cite{S_and_S} 
So we can safely conclude that for all $n>N_{85}$ the $k$th row of the \Sierpinski{} matrix $S_n$ certainly contains a prime  for $1 \leq k \leq \hbox{57 790}$.

\clearpage
\item 
Finally, another useful choice is to take  $m=3$, $\eta_3 = 0.15$, for which one has $x_3= \hbox{19 035 709 163} = p_{{{841\;508\;302}}}$ [Axler~\cite{Axler}].
The inequality is then valid only for $k\geq \ceil[x_3/N_{85}]=2$, and is  
easily checked to be satisfied on the interval $k\in[2,\; \hbox{141 618}]$. Of course the interval $k\in[1,9]$ was already dealt with in \Sierpinski's original discussion~\cite{S_and_S} so we can safely conclude that for all $n>N_{83}$ the $k$th row of the \Sierpinski{} matrix $S_n$ certainly contains a prime  for $1 \leq k \leq \hbox{141 618}$.
\end{itemize}
Consequently for all $n>N_{83}$ the $k$th row of the \Sierpinski{} matrix $S_n$ certainly contains a prime for $1 \leq k \leq \hbox{141 618}$.
\hfill $\Box$.

%================================================
\bigskip
\hrule\hrule\hrule
%================================================
\section{Discussion}
%================================================

We have seen that:
\begin{itemize}
\item 
\Sierpinski's Hypothesis H$_1$ holds at least up to  $n = \hbox{10 070 368 414} \gtrsim 10 \hbox{ billion}$. 
\item 
For arbitrary $n$ \Sierpinski's Hypothesis H$_1$ holds for at least 
one quarter of the rows of the \Sierpinski{} matrix.
\item
For arbitrary $n$ \Sierpinski's Hypothesis H$_1$ holds for at least  the first $\hbox{141 618}$ rows of the \Sierpinski{} matrix.
\end{itemize}
How could these results be improved?
\begin{itemize}
\item 
Most straightforwardly, locating and verifying a few more maximal prime gaps would certainly improve the range of validity.
\item 
More subtly, even developing some tractable bounds on the location and sizes of the maximal prime gaps  (even without demanding precise locations and sizes) would almost certainly improve the range of validity.
\item 
Any improvement on the Montgomery--Vaughn bound
$\pi(x+y) < \pi(x)+ {2y\over\ln y}$ would increase the fraction of rows of the \Sierpinski{} matrix guaranteed to contain primes.
\item
Improved bounds on the first Chebyshev function would almost certainly improve the state of knowledge regarding which of the lower rows of the \Sierpinski{} matrix could be guaranteed to contain primes.
\end{itemize}
Overall, while this current document certainly places some rather stringent bounds on  \Sierpinski's Hypothesis H$_1$, it is also clear that there are certainly many promising avenues for further development.

%%================================================
%\section*{Acknowledgement}
%%================================================

%=================================================
\bigskip
\hrule\hrule\hrule
\bigskip
%=================================================
%\bigskip
%\hrule\hrule\hrule
%================================================
%================================================
\clearpage
%================================================

%================================================
\end{document}